# TRIVIAL INTERSECTION OF $\sigma$-FIELDS AND GIBBS SAMPLING

By Patrizia Berti, Luca Pratelli and Pietro Rigo

*Università di Modena e Reggio-Emilia, Accademia Navale di Livorno and Università di Pavia*

Let $(\Omega, \mathcal{F}, P)$ be a probability space and $\mathcal{N}$ the class of those $F \in \mathcal{F}$ satisfying $P(F) \in \{0,1\}$. For each $\mathcal{G} \subset \mathcal{F}$, define $\overline{\mathcal{G}} = \sigma(\mathcal{G} \cup \mathcal{N})$. Necessary and sufficient conditions for $\overline{\mathcal{A}} \cap \overline{\mathcal{B}} = \overline{\mathcal{A} \cap \mathcal{B}}$, where $\mathcal{A}, \mathcal{B} \subset \mathcal{F}$ are sub-$\sigma$-fields, are given. These conditions are then applied to the (two-component) Gibbs sampler. Suppose $X$ and $Y$ are the coordinate projections on $(\Omega, \mathcal{F}) = (\mathcal{X} \times \mathcal{Y}, \mathcal{U} \otimes \mathcal{V})$ where $(\mathcal{X}, \mathcal{U})$ and $(\mathcal{Y}, \mathcal{V})$ are measurable spaces. Let $(X_n, Y_n)_{n \geq 0}$ be the Gibbs chain for $P$. Then, the SLLN holds for $(X_n, Y_n)$ if and only if $\overline{\sigma(X)} \cap \overline{\sigma(Y)} = \mathcal{N}$, or equivalently if and only if $P(X \in U)P(Y \in V) = 0$ whenever $U \in \mathcal{U}$, $V \in \mathcal{V}$ and $P(U \times V) = P(U^c \times V^c) = 0$. The latter condition is also equivalent to ergodicity of $(X_n, Y_n)$, on a certain subset $S_0 \subset \Omega$, in case $\mathcal{F} = \mathcal{U} \otimes \mathcal{V}$ is countably generated and $P$ absolutely continuous with respect to a product measure.

**1. The problem.** Let $(\Omega, \mathcal{F}, P)$ be a probability space and $\mathcal{A}, \mathcal{B} \subset \mathcal{F}$ sub-$\sigma$-fields. Letting

$$\mathcal{N} = \{F \in \mathcal{F} : P(F) \in \{0,1\}\} \quad \text{and} \quad \overline{\mathcal{G}} = \sigma(\mathcal{G} \cup \mathcal{N}),$$

for any subclass $\mathcal{G} \subset \mathcal{F}$, we aim at giving conditions for

(1) $$\overline{\mathcal{A}} \cap \overline{\mathcal{B}} = \overline{\mathcal{A} \cap \mathcal{B}}.$$

**2. Motivations.** Apart from its possible theoretical interest, there are three (nonindependent) reasons for investigating (1).

2.1. *Iterated conditional expectations.* Given a real random variable $Z$ satisfying $E\{|Z|\log(1+|Z|)\} < \infty$, define $Z_0 = Z$ and $\mathcal{G}_n = \mathcal{A}$ or $\mathcal{G}_n = \mathcal{B}$ as $n$ is even or odd. By a classical result of Burkholder and Chow [2] and Burkholder [4], one obtains

$$(2) \qquad Z_n := E(Z_{n-1} \,|\, \mathcal{G}_n) \to E(Z \,|\, \overline{\mathcal{A}} \cap \overline{\mathcal{B}}) \qquad \text{a.s.}$$

A natural question is whether $E(Z \,|\, \mathcal{A} \cap \mathcal{B})$ can be taken as the limit in (2).

COROLLARY 2.1. $Z_n \to E(Z \,|\, \mathcal{A} \cap \mathcal{B})$ a.s., for all real random variables $Z$ such that $E\{|Z|\log(1+|Z|)\} < \infty$, if and only if condition (1) holds.

PROOF. Under (1), just note that $E(Z \,|\, \mathcal{A} \cap \mathcal{B})$ is a version of $E(Z \,|\, \overline{\mathcal{A} \cap \mathcal{B}})$. Conversely, suppose $Z_n \to E(Z \,|\, \mathcal{A} \cap \mathcal{B})$ a.s. for all $Z$. Since $\overline{\mathcal{A} \cap \mathcal{B}} \subset \overline{\mathcal{A}} \cap \overline{\mathcal{B}}$, it suffices to prove that $\overline{\mathcal{A}} \cap \overline{\mathcal{B}} \subset \overline{\mathcal{A} \cap \mathcal{B}}$. Given $F \in \overline{\mathcal{A}} \cap \overline{\mathcal{B}}$, condition (2) implies $I_F = E(I_F \,|\, \overline{\mathcal{A}} \cap \overline{\mathcal{B}}) = E(I_F \,|\, \mathcal{A} \cap \mathcal{B})$ a.s. Letting $F_0 = \{E(I_F \,|\, \mathcal{A} \cap \mathcal{B}) \neq I_F\}$, and noting that $P(F_0) = 0$, yields

$$F = (F \cap F_0) \cup (\{E(I_F \,|\, \mathcal{A} \cap \mathcal{B}) = 1\} \cap F_0^c) \in \overline{\mathcal{A} \cap \mathcal{B}}. \qquad \square$$

As an application, think of a problem where $E(\cdot \,|\, \mathcal{A})$ and $E(\cdot \,|\, \mathcal{B})$ are easy to evaluate while $E(\cdot \,|\, \mathcal{A} \cap \mathcal{B})$ is not. In order to estimate $E(Z \,|\, \mathcal{A} \cap \mathcal{B})$, one strategy is using condition (2), but this is possible precisely when (1) holds.

2.2. *Sufficiency.* Suppose that, rather than a single probability measure $P$, we are given a collection $M$ of probability measures $Q$ on $(\Omega, \mathcal{F})$. For any $\mathcal{G} \subset \mathcal{F}$ define $\mathcal{G}_M = \sigma(\mathcal{G} \cup \mathcal{N}_M)$, where $\mathcal{N}_M$ is the class of those $F \in \mathcal{F}$ such that $Q(F) = 0$ for all $Q \in M$ or $Q(F) = 1$ for all $Q \in M$. In this framework, condition (1) turns into

$$(1^*) \qquad \mathcal{A}_M \cap \mathcal{B}_M = (\mathcal{A} \cap \mathcal{B})_M.$$

A sub-$\sigma$-field $\mathcal{G} \subset \mathcal{F}$ is *sufficient* (for $M$) in case, for each $F \in \mathcal{F}$, there is a $\mathcal{G}$-measurable function $f : \Omega \to \mathbb{R}$ which is a version of $E_Q(I_F \,|\, \mathcal{G})$ for all $Q \in M$.

Generally, sufficiency of both $\mathcal{A}$ and $\mathcal{B}$ does not imply that of $\mathcal{A} \cap \mathcal{B}$. By Theorem 4 of [3], however, $\mathcal{A} \cap \mathcal{B}$ is sufficient provided $\mathcal{A}$ and $\mathcal{B}$ are sufficient and at least one of them includes $\mathcal{N}_M$. This implies the following result.

COROLLARY 2.2. $\mathcal{A} \cap \mathcal{B}$ *is sufficient whenever* $\mathcal{A}$ *and* $\mathcal{B}$ *are sufficient and condition* $(1^*)$ *holds.*

PROOF. We first verify that $\mathcal{G}$ is sufficient if and only if $\mathcal{G}_M$ is sufficient, where $\mathcal{G} \subset \mathcal{F}$ is any sub-$\sigma$-field. The "only if" part is trivial. Suppose $\mathcal{G}_M$ is sufficient, fix $F \in \mathcal{F}$, and take a $\mathcal{G}_M$-measurable function $f$ which is a version of $E_Q(I_F \,|\, \mathcal{G}_M)$ for all $Q \in M$. Since

$$\mathcal{G}_M = \{F \in \mathcal{F} : \text{there is } G \in \mathcal{G} \text{ such that } Q(F \Delta G) = 0 \text{ for all } Q \in M\},$$

for each $n$ there is a $\mathcal{G}$-measurable function $\phi_n$ such that $Q(|f - \phi_n| < \frac{1}{n}) = 1$ for all $Q \in M$. For $\omega \in \Omega$, define $\phi(\omega) = \lim_n \phi_n(\omega)$ if the limit exists and $\phi(\omega) = 0$ otherwise. Then, $\phi$ is $\mathcal{G}$-measurable and $Q(f = \phi) = 1$ for all $Q \in M$. Thus, $\phi$ is a version of $E_Q(I_F \mid \mathcal{G})$ for all $Q \in M$, which shows that $\mathcal{G}$ is sufficient. Next, since $\mathcal{A}$ and $\mathcal{B}$ are sufficient, $\mathcal{A}_M$ and $\mathcal{B}_M$ are still sufficient, and thus $(\mathcal{A} \cap \mathcal{B})_M = \mathcal{A}_M \cap \mathcal{B}_M$ is sufficient by Theorem 4 of [3]. Therefore, $\mathcal{A} \cap \mathcal{B}$ is sufficient. $\square$

2.3. *Two-component Gibbs sampler.* Suppose $(\Omega, \mathcal{F}) = (\mathcal{X} \times \mathcal{Y}, \mathcal{U} \otimes \mathcal{V})$ is the product of two measurable spaces $(\mathcal{X}, \mathcal{U})$ and $(\mathcal{Y}, \mathcal{V})$ and let $X : \Omega \to \mathcal{X}$, $Y : \Omega \to \mathcal{Y}$ be the coordinate projections. Suppose also that regular versions of the conditional distribution of $Y$ given $X$ and $X$ given $Y$ are available under $P$ (precise definitions are given in Section 4). Roughly speaking, the Gibbs-chain $(X_n, Y_n)_{n \geq 0}$ can be described as follows. Starting from $\omega = (x, y)$, the next state $\omega^* = (x^*, y^*)$ is obtained by first choosing $y^*$ from the conditional distribution of $Y$ given $X = x$ and then $x^*$ from the conditional distribution of $X$ given $Y = y^*$. Iterating this procedure produces a homogeneous Markov chain $(X_n, Y_n)$ with stationary distribution $P$. Let $\mathbb{P}$ denote the law of such a chain when $(X_0, Y_0) \sim P$, and let

$$m_n(\phi) = \frac{1}{n} \sum_{i=0}^{n-1} \phi(X_i, Y_i)$$

for each function $\phi : \Omega \to \mathbb{R}$.

In real problems, $(X_n, Y_n)$ is constructed mainly for sampling from $P$. To this end, it is crucial that the SLLN is available under $\mathbb{P}$, that is,

(3) $$m_n(\phi) \to \int \phi \, dP \quad \mathbb{P}\text{-a.s., for all } \phi \in L_1(P).$$

Note that, under (3), for each probability measure $Q \ll P$ one also obtains

$$m_n(\phi) \to \int \phi \, dP \quad \mathbb{Q}\text{-a.s., for each } \phi \in L_1(P)$$

where $\mathbb{Q}$ is the law of the chain $(X_n, Y_n)$ when $(X_0, Y_0) \sim Q$.

In addition to (3), various other properties are usually requested to $(X_n, Y_n)$, mainly ergodicity, CLT and the convergence rate. Nevertheless, condition (3) seems (to us) a fundamental one. It is a sort of necessary condition, since the Gibbs sampling procedure does not make sense without (3). Accordingly, we say that $P$ is *Gibbs-admissible* in case (3) holds.

But what about condition (1)? The link is that $P$ turns out to be Gibbs-admissible precisely when

$$\overline{\sigma(X)} \cap \overline{\sigma(Y)} = \mathcal{N} = \overline{\sigma(X) \cap \sigma(Y)}.$$

In other terms, the Gibbs sampling procedure makes sense for $P$ (in the SLLN-sense) if and only if $P$ meets condition (1) with $\mathcal{A} = \sigma(X)$ and $\mathcal{B} = \sigma(Y)$.

In fact, something more is true. Let $K$ be the transition kernel of $(X_n, Y_n)$ and $S_0 = \{\omega \in \Omega : K(\omega, \cdot) \ll P\}$. Under mild conditions ($\mathcal{F}$ countably generated and $P$ absolutely continuous with respect to a product measure), one obtains $P(S_0) = 1$ and

$$\overline{\sigma(X)} \cap \overline{\sigma(Y)} = \mathcal{N} \quad \Leftrightarrow \quad (X_n, Y_n) \text{ is ergodic on } S_0.$$

For proving the previous statements, a key ingredient is a result of Diaconis et al., connecting the Gibbs-chain $(X_n, Y_n)$ with the Burkholder–Chow result of Section 2.1; see Theorem 4.1 of [6]. Indeed, $\overline{\sigma(X)} \cap \overline{\sigma(Y)} = \mathcal{N}$ appears as an assumption in various results from [6] (Corollary 3.1, Propositions 5.1, 5.2 and 5.3). Also, investigating $\overline{\sigma(X)} \cap \overline{\sigma(Y)} = \mathcal{N}$ was suggested to us by Persi Diaconis.

One of our main results (Corollary 3.5) is that $\overline{\sigma(X)} \cap \overline{\sigma(Y)} = \mathcal{N}$ is equivalent to the following property of $P$:

$$P(X \in U) = 0 \quad \text{or} \quad P(Y \in V) = 0$$
$$\text{whenever } U \in \mathcal{U}, V \in \mathcal{V} \text{ and } P(U \times V) = P(U^c \times V^c) = 0.$$

This paper is organized into two parts. Section 3 gives general results on conditions (1) and (1\*). It includes characterizations, examples and various working sufficient conditions in case $P$ is absolutely continuous with respect to a product measure. Under this assumption, it is also shown that $P$ is *atomic on* $\overline{\sigma(X)} \cap \overline{\sigma(Y)}$ (with $X$ and $Y$ as in Section 2.3). The main results are Theorems 3.1 and 3.10 and Corollaries 3.5, 3.7 and 3.9. Section 4 deals with the Gibbs sampler and contains the material sketched above. The main results are Theorems 4.2, 4.4 and 4.5.

**3. General results.** This section is split into three subsections. All examples are postponed to the last one.

3.1. *Necessary and sufficient conditions.* Condition (1) admits a surprisingly simple characterization.

THEOREM 3.1. *Let*

$$\mathcal{J} = \{A \cap B : A \in \mathcal{A}, B \in \mathcal{B} \text{ and } P(A \cap B) + P(A^c \cap B^c) = 1\}.$$

*Then,* $\overline{\mathcal{A} \cap \mathcal{B}} = \overline{\mathcal{J}}$. *Moreover,* $\overline{\mathcal{A}} \cap \overline{\mathcal{B}} = \overline{\mathcal{A} \cap \mathcal{B}}$ *if and only if*

(4)
$$A \in \mathcal{A}, B \in \mathcal{B} \quad \text{and} \quad P(A \cap B) = P(A^c \cap B^c) = 0$$
$$\text{implies } P(A \Delta D) = 0 \text{ or } P(B \Delta D) = 0 \text{ for some } D \in \mathcal{A} \cap \mathcal{B}.$$

PROOF. First note that, for any sub-$\sigma$-field $\mathcal{G} \subset \mathcal{F}$, one has
$$\overline{\mathcal{G}} = \{F \in \mathcal{F} : P(F \Delta G) = 0 \text{ for some } G \in \mathcal{G}\}.$$
Let $F \in \overline{\mathcal{A}} \cap \overline{\mathcal{B}}$. Then, $P(A \Delta F) = P(B \Delta F) = 0$, for some $A \in \mathcal{A}$ and $B \in \mathcal{B}$, and
$$1 - P(A \cap B) - P(A^c \cap B^c) = P(A \Delta B) \leq P(A \Delta F) + P(B \Delta F) = 0.$$
Hence, $J := A \cap B \in \mathcal{J}$. Since $P(F \Delta J) \leq P(A \Delta F) + P(B \Delta F) = 0$, then $F \in \overline{\mathcal{J}}$. Conversely, let $J = A \cap B \in \mathcal{J}$ where $A \in \mathcal{A}$, $B \in \mathcal{B}$ and $P(A \cap B) + P(A^c \cap B^c) = 1$. Define $H = (A \cap B) \cup (A^c \cap B^c)$. Since $P(H) = 1$ and $J = A \cap H = B \cap H$, it follows that $J \in \overline{\mathcal{A}} \cap \overline{\mathcal{B}}$. Therefore, $\overline{\mathcal{A}} \cap \overline{\mathcal{B}} = \overline{\mathcal{J}}$. In particular, $\overline{\mathcal{A}} \cap \overline{\mathcal{B}} = \overline{\mathcal{A} \cap \mathcal{B}}$ if and only if $\mathcal{J} \subset \overline{\mathcal{A} \cap \mathcal{B}}$, and thus it remains only to show that $\mathcal{J} \subset \overline{\mathcal{A} \cap \mathcal{B}}$ is equivalent to condition (4). Suppose (4) holds and fix $J \in \mathcal{J}$. Then, $J$ can be written as $J = A \cap B^c$ for some $A \in \mathcal{A}$ and $B \in \mathcal{B}$ with $P(A \cap B) = P(A^c \cap B^c) = 0$. By (4), it follows that $A \in \overline{\mathcal{A} \cap \mathcal{B}}$ or $B \in \overline{\mathcal{A} \cap \mathcal{B}}$, say $A \in \overline{\mathcal{A} \cap \mathcal{B}}$. Since $P(A \cap B) = 0$, one obtains $J = A \cap B^c = A - (A \cap B) \in \overline{\mathcal{A} \cap \mathcal{B}}$. Finally, suppose $\mathcal{J} \subset \overline{\mathcal{A} \cap \mathcal{B}}$ and fix $A \in \mathcal{A}$ and $B \in \mathcal{B}$ with $P(A \cap B) = P(A^c \cap B^c) = 0$. Since $A \cap B^c \in \mathcal{J} \subset \overline{\mathcal{A} \cap \mathcal{B}}$ and $P(A \cap B) = 0$, then $A = (A \cap B) \cup (A \cap B^c) \in \overline{\mathcal{A} \cap \mathcal{B}}$, that is, $P(A \Delta D) = 0$ for some $D \in \mathcal{A} \cap \mathcal{B}$. $\square$

The case of $k \geq 2$ sub-$\sigma$-fields can be settled essentially as in Theorem 3.1.

THEOREM 3.2. *Let $\mathcal{A}_1, \ldots, \mathcal{A}_k \subset \mathcal{F}$ be sub-$\sigma$-fields, $k \geq 2$. Then,*
$$\bigcap_{i=1}^{k} \overline{\mathcal{A}_i} = \overline{\bigcap_{i=1}^{k} \mathcal{A}_i}$$
*if and only if*
$$A_1 \in \mathcal{A}_1, \ldots, A_k \in \mathcal{A}_k \quad \text{and} \quad P\left(\bigcap_{i=1}^{k} A_i\right) + P\left(\bigcap_{i=1}^{k} A_i^c\right) = 1$$
*implies $P(A_i \Delta D) = 0$ for some $i$ and $D \in \bigcap_{i=1}^{k} \mathcal{A}_i$.*

The simple argument which leads to Theorem 3.1 allows to find conditions for sufficiency of $\mathcal{A} \cap \mathcal{B}$ as well.

THEOREM 3.3. *In the notation of Section 2.2, $\mathcal{A}_M \cap \mathcal{B}_M = (\mathcal{A} \cap \mathcal{B})_M$ if and only if*

$$A \in \mathcal{A}, B \in \mathcal{B} \quad \text{and} \quad Q(A \cap B) = Q(A^c \cap B^c) = 0 \quad \text{for all } Q \in M$$

(4*) *implies the existence of $D \in \mathcal{A} \cap \mathcal{B}$ such that*

$$Q(A \Delta D) = 0 \quad \text{for all } Q \in M \quad \text{or} \quad Q(B \Delta D) = 0 \quad \text{for all } Q \in M.$$

Hence, by Corollary 2.2 and Theorem 3.3, $\mathcal{A} \cap \mathcal{B}$ is sufficient whenever $\mathcal{A}$ and $\mathcal{B}$ are sufficient and condition (4*) holds.

The proofs of both Theorems 3.2 and 3.3 have been omitted since they are quite analogous to that of Theorem 3.1.

Let us come back to the main concern of this paper, that is, a single probability measure $P$ and two sub-$\sigma$-fields $\mathcal{A}$ and $\mathcal{B}$. Condition (4) of Theorem 3.1 trivially holds provided

$$A \in \mathcal{A}, B \in \mathcal{B} \quad \text{and} \quad P(A \cap B) = P(A^c \cap B^c) = 0$$
(5)
$$\text{implies } P(A) = 0 \text{ or } P(B) = 0.$$

Generally, condition (5) is stronger than (4) [just take $\mathcal{A} = \mathcal{B}$, so that (4) trivially holds, and choose $P$ such that (5) fails]. However, (4) and (5) are equivalent in a significant particular case.

COROLLARY 3.4. *If $\mathcal{A} \cap \mathcal{B} \subset \mathcal{N}$, then*

$$\overline{\mathcal{A}} \cap \overline{\mathcal{B}} = \mathcal{N} \quad \Leftrightarrow \quad \text{condition (4) holds} \quad \Leftrightarrow \quad \text{condition (5) holds}.$$

PROOF. Suppose (4) holds and fix $A \in \mathcal{A}$, $B \in \mathcal{B}$ with $P(A \cap B) = P(A^c \cap B^c) = 0$. Then, $P(A \Delta D) P(B \Delta D) = 0$ for some $D \in \mathcal{A} \cap \mathcal{B}$, say $P(A \Delta D) = 0$. Since $\mathcal{A} \cap \mathcal{B} \subset \mathcal{N}$, one obtains $P(A) = P(D) \in \{0, 1\}$. If $P(A) = 0$, then $P(A)P(B) = 0$. If $P(A) = 1$, then $P(B) = P(A \cap B) = 0$ and again $P(A)P(B) = 0$. Thus, (5) holds. Since $\overline{\mathcal{A} \cap \mathcal{B}} = \mathcal{N}$, an application of Theorem 3.1 concludes the proof. □

In the sequel, we deal with product measurable spaces. Let $(\mathcal{X}, \mathcal{U})$ and $(\mathcal{Y}, \mathcal{V})$ be measurable spaces and

$$(\Omega, \mathcal{F}) = (\mathcal{X} \times \mathcal{Y}, \mathcal{U} \otimes \mathcal{V}), \qquad X(x, y) = x, \qquad Y(x, y) = y,$$

where $(x, y) \in \mathcal{X} \times \mathcal{Y}$. We focus on $\mathcal{A} = \sigma(X)$ and $\mathcal{B} = \sigma(Y)$ and we let

$$\mathcal{D} = \overline{\sigma(X)} \cap \overline{\sigma(Y)}.$$

On noting that $\sigma(X) \cap \sigma(Y) = \{\varnothing, \Omega\}$, Corollary 3.4 implies the following statement.

COROLLARY 3.5. $\mathcal{D} = \mathcal{N}$ *if and only if*

$$U \in \mathcal{U}, V \in \mathcal{V} \quad \text{and} \quad P(U \times V) = P(U^c \times V^c) = 0$$
(6)
$$\textit{implies } P(X \in U) = 0 \textit{ or } P(Y \in V) = 0.$$

By Corollary 3.5, if $\mathcal{D} \neq \mathcal{N}$, then $\mathcal{D}$ includes a rectangle $U \times V$ such that $U \in \mathcal{U}$, $V \in \mathcal{V}$ and $0 < P(U \times V) < 1$. This fact implies a first sufficient condition for $\mathcal{D} = \mathcal{N}$.

COROLLARY 3.6. *For $\mathcal{D} = \mathcal{N}$, it is sufficient that*

(7)
$$E(E(I_{U \times V} \mid \sigma(X)) \mid \sigma(Y)) + E(E(I_{U \times V} \mid \sigma(Y)) \mid \sigma(X)) > 0 \quad \text{a.s.}$$
$$\text{whenever } U \in \mathcal{U}, V \in \mathcal{V} \text{ and } 0 < P(U \times V) < 1.$$

PROOF. Let $U \in \mathcal{U}$ and $V \in \mathcal{V}$. If $U \times V \in \mathcal{D}$, then

$$E(E(I_{U \times V} \mid \sigma(X)) \mid \sigma(Y)) + E(E(I_{U \times V} \mid \sigma(Y)) \mid \sigma(X)) = 2 I_{U \times V} \quad \text{a.s.}$$

Thus, (7) implies $U \times V \notin \mathcal{D}$ in case $0 < P(U \times V) < 1$. $\square$

3.2. *Sufficient conditions in case $P$ is absolutely continuous with respect to a product measure.* In real problems, $P$ usually has a density with respect to some product measure. Let $\mu$ and $\nu$ be $\sigma$-finite measures on $\mathcal{U}$ and $\mathcal{V}$, respectively. In this subsection, $P \ll \mu \times \nu$ and $f$ is a version of the density of $P$ with respect to $\mu \times \nu$.

COROLLARY 3.7. *For $\mathcal{D} = \mathcal{N}$, it is sufficient that $P \ll \mu \times \nu$ and*

$$(U_0 \times \mathcal{Y}) \cup (\mathcal{X} \times V_0) \supset \{f > 0\} \supset U_0 \times V_0$$

*for some $U_0 \in \mathcal{U}$ and $V_0 \in \mathcal{V}$ such that $P(U_0 \times V_0) > 0$.*

PROOF. By Corollary 3.5, it suffices to prove condition (6). Let $U \in \mathcal{U}$ and $V \in \mathcal{V}$ be such that $P(U \times V) = P(U^c \times V^c) = 0$. Since $f > 0$ on $U_0 \times V_0$ and

$$\int_{U \cap U_0 \times V \cap V_0} f \, d(\mu \times \nu) = P((U \cap U_0) \times (V \cap V_0)) \leq P(U \times V) = 0,$$

it follows that

$$\mu(U \cap U_0)\nu(V \cap V_0) = \mu \times \nu((U \cap U_0) \times (V \cap V_0)) = 0,$$

say $\mu(U \cap U_0) = 0$. Similarly, $P(U^c \times V^c) = 0$ and $f > 0$ on $U_0 \times V_0$ imply $\mu(U_0 - U)\nu(V_0 - V) = 0$. Since $P(X \in U_0) > 0$ and $\mu(U \cap U_0) = 0$, it must be $\mu(U_0 - U) > 0$ and thus $\nu(V_0 - V) = 0$. Let $H_0 = (U_0 \times \mathcal{Y}) \cup (\mathcal{X} \times V_0)$. Since $P(H_0) = 1$ and $P(X \in U \cap U_0) = P(Y \in V_0 - V) = 0$, one obtains

$$P(X \in U) = P(\{X \in U - U_0\} \cap H_0) = P((U - U_0) \times V_0)$$
$$= P((U - U_0) \times (V \cap V_0)) \leq P(U \times V) = 0.$$

Therefore, condition (6) holds. $\square$

Corollary 3.7 applies in particular if $\{f > 0\} \supset U_0 \times \mathcal{Y}$ for some $U_0 \in \mathcal{U}$ with $P(X \in U_0) > 0$ (just take $V_0 = \mathcal{Y}$). Likewise, it applies if $\{f > 0\} \supset \mathcal{X} \times V_0$ for some $V_0 \in \mathcal{V}$ such that $P(Y \in V_0) > 0$.

Let $\mu_0$ be a probability measure on $\mathcal{U}$ equivalent to $\mu$, that is, $\mu_0 \ll \mu$ and $\mu \ll \mu_0$. Similarly, let $\nu_0$ be a probability measure on $\mathcal{V}$ equivalent to $\nu$. Then, $\mu_0 \times \nu_0$ is equivalent to $\mu \times \nu$ and, for each $H \in \mathcal{F}$ with $\mu \times \nu(H) > 0$, one can define the probability measure

$$Q_H(F) = \frac{\mu_0 \times \nu_0(F \cap H)}{\mu_0 \times \nu_0(H)}, \qquad F \in \mathcal{F}.$$

Say that $H$ has the *trivial intersection property*, or briefly that $H$ is TIP, in case $H \in \mathcal{F}$, $\mu \times \nu(H) > 0$ and $\mathcal{D} = \mathcal{N}$ holds when $P = Q_H$. Note that whether or not $H$ is TIP does not depend on the choice of $\mu_0$ and $\nu_0$. Note also that

$$\mathcal{D} = \mathcal{N} \quad \Leftrightarrow \quad \text{the set } \{f > 0\} \text{ is TIP}.$$

COROLLARY 3.8. *Suppose $P \ll \mu \times \nu$ and $\{f > 0\} = \bigcup_n H_n$, where $H_1 \subset H_2 \subset \cdots$ is an increasing sequence of TIP sets. Then, $\mathcal{D} = \mathcal{N}$.*

PROOF. Let $H = \{f > 0\}$ and $U \in \mathcal{U}, V \in \mathcal{V}$ with $Q_H(U \times V) = Q_H(U^c \times V^c) = 0$. Since $H_n \subset H$, then $Q_{H_n}(U \times V) = Q_{H_n}(U^c \times V^c) = 0$ for all $n$. Since $H_n$ is TIP, Corollary 3.5 (applied to $Q_{H_n}$) yields $Q_{H_n}(X \in U)Q_{H_n}(Y \in V) = 0$ for all $n$. Thus, $H_n \uparrow H$ implies

$$Q_H(X \in U)Q_H(Y \in V) = \lim_n Q_{H_n}(X \in U)Q_{H_n}(Y \in V) = 0.$$

By Corollary 3.5, $H = \{f > 0\}$ is TIP. □

In order to generalize Corollary 3.8, one more piece of terminology is useful. Given two sets $F, G \in \mathcal{F}$, say that $F$ *communicates with* $G$ in case at least one of the following conditions (i) and (ii) holds:

(i) there is $V_0 \in \mathcal{V}$ with $\nu(V_0) > 0$ and $\mu(F^y) > 0$, $\mu(G^y) > 0$ for all $y \in V_0$,
(ii) there is $U_0 \in \mathcal{U}$ with $\mu(U_0) > 0$ and $\nu(F_x) > 0$, $\nu(G_x) > 0$ for all $x \in U_0$,

where $F_x = \{y \in \mathcal{Y} : (x, y) \in F\}$ and $F^y = \{x \in \mathcal{X} : (x, y) \in F\}$ denote the sections of $F$.

COROLLARY 3.9. *Suppose $P \ll \mu \times \nu$ and $\{f > 0\} = \bigcup_n H_n$, where $H_n$ is TIP and $H_n$ communicates with $H_{n+1}$ for each $n$. Then, $\mathcal{D} = \mathcal{N}$.*

PROOF. It is enough to prove that $F \cup G$ is TIP whenever $F$ and $G$ are TIP and $F$ communicates with $G$. In that case, in fact, since $\bigcup_{i=1}^{n-1} H_i$ communicates with $H_n$, a simple induction implies that $\bigcup_{i=1}^{n} H_i$ is TIP for all $n$. Hence, $\mathcal{D} = \mathcal{N}$ by Corollary 3.8.

Suppose $F$ and $G$ are TIP and condition (i) holds [the proof is the same if (ii) holds]. Set $H = F \cup G$ and fix $U \in \mathcal{U}$, $V \in \mathcal{V}$ with $Q_H(U \times V) = Q_H(U^c \times V^c) = 0$. Since $F$ and $G$ are TIP and

$$Q_F(U \times V) = Q_F(U^c \times V^c) = Q_G(U \times V) = Q_G(U^c \times V^c) = 0,$$

one obtains

$$Q_F(X \in U) = 0 \text{ or } Q_F(Y \in V) = 0 \quad \text{and}$$
$$Q_G(X \in U) = 0 \text{ or } Q_G(Y \in V) = 0.$$

Let $V_0$ be as in condition (i). If $Q_F(X \in U) = 0$, then $Q_F(Y \in V) = 1$. By (i) and since $\mu_0$ and $\nu_0$ are equivalent to $\mu$ and $\nu$, it follows that

$$Q_F(Y \in V \cap V_0) = Q_F(Y \in V_0) = \frac{\int_{V_0} \mu_0(F^y) \nu_0(dy)}{\mu_0 \times \nu_0(F)} > 0.$$

Hence $\nu_0(V \cap V_0) > 0$, and this implies

$$Q_G(Y \in V) \geq Q_G(Y \in V \cap V_0) = \frac{\int_{V \cap V_0} \mu_0(G^y) \nu_0(dy)}{\mu_0 \times \nu_0(G)} > 0.$$

Therefore, $Q_F(X \in U) = 0$ implies $Q_G(Y \in V) > 0$, and similarly $Q_G(X \in U) = 0$ implies $Q_F(Y \in V) > 0$. It follows that $Q_F(X \in U) = Q_G(X \in U) = 0$ or $Q_F(Y \in V) = Q_G(Y \in V) = 0$, which implies $Q_H(X \in U) = 0$ or $Q_H(Y \in V) = 0$. Thus, condition (6) holds for $Q_H$, and $H = F \cup G$ is TIP by Corollary 3.5. $\square$

So far, conditions for $P$ to be 0–1-valued on $\mathcal{D}$ have been given. A weaker but useful result is that the latter property holds locally, in the sense that $\Omega$ can be partitioned into sets $H_1, H_2, \ldots \in \mathcal{D}$ such that $P(H_n) > 0$ and $P(\cdot \mid H_n)$ is 0–1-valued on $\mathcal{D}$ for each $n$. We now prove that this is always true provided $P \ll \mu \times \nu$. In that case, in fact, $P$ is *atomic* on $\mathcal{D}$. Recall that, given a probability space $(\mathcal{Z}, \mathcal{E}, Q)$, a $Q$-atom is a set $K \in \mathcal{E}$ such that $Q(K) > 0$ and $Q(\cdot \mid K)$ is 0–1-valued. In general, there are three possible situations: (j) $Q$ is nonatomic, that is, there are no $Q$-atoms; (jj) $Q$ is atomic, that is, the $Q$-atoms form a partition of $\mathcal{Z}$; (jjj) there is $K \in \mathcal{E}$, $0 < Q(K) < 1$, such that $Q(\cdot \mid K)$ is nonatomic and $K^c$ is a disjoint union of $Q$-atoms.

THEOREM 3.10. *If $P \ll \mu \times \nu$, then $P$ is atomic on $\mathcal{D}$ (i.e., the restriction of $P$ to $\mathcal{D}$ is atomic).*

PROOF. Fix $H \in \mathcal{D}$ with $P(H) > 0$ and let $P_H$ denote the restriction of $P(\cdot \mid H)$ to $\mathcal{D}$. If $P_H$ is nonatomic, the probability space $(\Omega, \mathcal{D}, P_H)$ supports a real random variable with uniform distribution on $(0, 1)$; see, for example, Theorem 3.1 of [1]. Hence, it suffices to prove that each $\mathcal{D}$-measurable function $Z:\Omega \to \mathbb{R}$ satisfies $P_H(Z \in C) = 1$ for some countable set $C \subset \mathbb{R}$. Let $Z:\Omega \to \mathbb{R}$ be $\mathcal{D}$-measurable. Since $\sigma(Z) \subset \overline{\sigma(X)}$, there is a $\mathcal{U}$-measurable function $h:\mathcal{X} \to \mathbb{R}$ satisfying $Z = h(X)$ a.s. Similarly, $\sigma(Z) \subset \overline{\sigma(Y)}$ yields $Z = k(Y)$ a.s. for some $\mathcal{V}$-measurable function $k:\mathcal{Y} \to \mathbb{R}$. Let $C = \{c \in \mathbb{R} : \nu\{y : k(y) = c\} > 0\}$. Since $\nu$ is $\sigma$-finite, $C$ is countable and

$$\mu \times \nu(h(X) \notin C, h(X) = k(Y)) = \int_{\{x \,:\, h(x) \notin C\}} \nu\{y : k(y) = h(x)\} \mu(dx) = 0.$$

Since $P \ll \mu \times \nu$ and $Z = h(X) = k(Y)$ a.s., it follows that

$$P(Z \in C) = 1 - P(h(X) \notin C, h(X) = k(Y)) = 1.$$

Thus, $P_H \ll P$ implies $P_H(Z \in C) = 1$. This concludes the proof. □

3.3. *Examples.* In this subsection, $\mathcal{X}$ and $\mathcal{Y}$ are topological spaces and $\mathcal{U}$ and $\mathcal{V}$ the corresponding Borel $\sigma$-fields. Moreover, $\mu$ and $\nu$ have full topological support (i.e., they are strictly positive on nonempty open sets) and $P$ has a density $f$ with respect to $\mu \times \nu$.

We note that, since $\mu$ and $\nu$ have full topological support, $F$ communicates with $G$ whenever $F, G \in \mathcal{F}$ and $F \cap G$ has nonempty interior. Further, by Corollary 3.9 (see also its proof), $F \cup G$ is TIP whenever $F$ and $G$ are TIP and $F$ communicates with $G$.

EXAMPLE 3.11. Let $\mathcal{X}$ and $\mathcal{Y}$ be second countable topological spaces. If $\{f > 0\}$ *is open and connected, then* $\mathcal{D} = \mathcal{N}$.

Suppose in fact that $H \subset \Omega$ is open and connected. Since $H$ is open, $H = \bigcup_n H_n$ where each $H_n$ is open and TIP (for instance, take the $H_n$ as open rectangles). For $\omega_1, \omega_2 \in H$, say that $\omega_1 \sim \omega_2$ in case there are a finite number of indices $j_1, \ldots, j_n$ such that $\omega_1 \in H_{j_1}$, $\omega_2 \in H_{j_n}$ and $H_{j_i} \cap H_{j_{i+1}} \neq \varnothing$ for each $i$. Then, $\sim$ is an equivalence relation on $H$. Since $H$ is connected and the equivalence classes of $\sim$ are open, there is precisely one equivalence class, that is, $\omega_1 \sim \omega_2$ for all $\omega_1, \omega_2 \in H$. Fix $\omega_0 \in H$. For each $k$, take $\omega_k \in H_k$ and define $M_k = H_k \cup (\bigcup_{i=1}^n H_{j_i})$, where $j_1, \ldots, j_n$ are such that $\omega_k \in H_{j_1}$, $\omega_0 \in H_{j_n}$ and $H_{j_i} \cap H_{j_{i+1}} \neq \varnothing$ for all $i$. By Corollary 3.9, $M_k$ is TIP. Thus, $H = \bigcup_k M_k$ is TIP as well, since $\omega_0 \in M_k \cap M_{k+1}$ for each $k$.

Note that $\{f > 0\}$ is open provided $f$ is lower semicontinuous. Thus, when $f$ is lower semicontinuous (and $\mathcal{X}$, $\mathcal{Y}$ are second countable and locally connected), a sufficient condition for $\mathcal{D} = \mathcal{N}$ is that the connected components $H_1, H_2, \ldots$ of $\{f > 0\}$ can be arranged in such a way that $H_n$ communicates with $H_{n+1}$ for all $n$. This follows from Example 3.11 and Corollary 3.9. Note also that, in case $\mathcal{X} = \mathbb{R}^n$, $\mathcal{Y} = \mathbb{R}^m$ and $\mu$, $\nu$ the Lebesgue measures, $P$ admits a lower semicontinuous density as far as it admits a Riemann-integrable density.

EXAMPLE 3.12. Let $\mathcal{X} = \mathbb{R}^n$, $\mathcal{Y} = \mathbb{R}^m$ and $\mu$, $\nu$ the Lebesgue measures. If $\{f > 0\}$ *is convex, then* $\mathcal{D} = \mathcal{N}$.

Suppose in fact that $H \in \mathcal{F}$ is convex and $\mu \times \nu(H) > 0$. Since $H$ is convex, $\mu \times \nu(H - H^0) \leq \mu \times \nu(\partial H) = 0$ where $H^0$ and $\partial H$ are the interior and the boundary of $H$. Thus, it suffices to note that $H^0$ is open and connected.

In various frameworks, for instance in Bayesian statistics, $P$ is usually a *mixture* of some other probability laws.

EXAMPLE 3.13. Let $(\Theta, \mathcal{E}, \pi)$ be a probability space and $\{P_\theta : \theta \in \Theta\}$ a collection of probabilities on $\mathcal{F}$ such that $\theta \mapsto P_\theta(H)$ is $\mathcal{E}$-measurable for fixed $H \in \mathcal{F}$. Define

$$P = \int P_\theta \pi(d\theta).$$

Then, condition (6) can be written as

(8) $\quad U \in \mathcal{U}, V \in \mathcal{V} \quad \text{and} \quad P_\theta(U \times V) = P_\theta(U^c \times V^c) = 0, \qquad \pi\text{-a.s.},$

$\quad\quad\quad$ implies $P_\theta(X \in U) = 0$, $\pi$-a.s. or $P_\theta(Y \in V) = 0$, $\pi$-a.s.

Several sufficient conditions for $\mathcal{D} = \mathcal{N}$ can be deduced from (8). For instance, $\mathcal{D} = \mathcal{N}$ provided each $P_\theta$ meets condition (6) [i.e., (6) holds when $P = P_\theta$] and

(9) $\quad \pi\{\theta : P_\theta(U \times V) = 1\} > 0 \quad \Longrightarrow \quad P_\theta(U \times V) > 0, \qquad \pi\text{-a.s.},$

for all $U \in \mathcal{U}$ and $V \in \mathcal{V}$. Suppose in fact $P_\theta(U \times V) = P_\theta(U^c \times V^c) = 0$, $\pi$-a.s., for some $U \in \mathcal{U}$ and $V \in \mathcal{V}$. Since each $P_\theta$ meets (6), then $P_\theta(U \times V^c) \in \{0, 1\}$, $\pi$-a.s. If $P_\theta(U \times V^c) = 0$, $\pi$-a.s., then $P_\theta(X \in U) = 0$, $\pi$-a.s. Otherwise, if

$$\pi\{\theta : P_\theta(U \times V^c) = 1\} = \pi\{\theta : P_\theta(U \times V^c) > 0\} > 0,$$

condition (9) yields $P_\theta(U \times V^c) > 0$, $\pi$-a.s. Hence, $P_\theta(Y \in V) = 0$, $\pi$-a.s. One more sufficient condition for $\mathcal{D} = \mathcal{N}$ applies when each $P_\theta$ has a density $f_\theta$ with respect to $\mu \times \nu$. In that case, $\mathcal{D} = \mathcal{N}$ whenever

(10) $\quad \{f_\theta > 0\}$ is TIP $\quad$ for all $\theta$ $\quad$ and

$\quad\quad\quad \{f_{\theta_1} > 0\}$ communicates with $\{f_{\theta_2} > 0\} \quad$ for all $\theta_1, \theta_2$.

Suppose in fact $P_{\theta_i}(U \times V) = P_{\theta_i}(U^c \times V^c) = 0$, $i = 1, 2$, for some $U \in \mathcal{U}$, $V \in \mathcal{V}$ and $\theta_1, \theta_2 \in \Theta$. By (10),

$$\{f_{\theta_1} > 0\} \cup \{f_{\theta_2} > 0\} \text{ is TIP.}$$

Using this fact, it is straightforward to verify that $P_{\theta_1}(X \in U) = P_{\theta_2}(X \in U) = 0$ or $P_{\theta_1}(Y \in V) = P_{\theta_2}(Y \in V) = 0$. Therefore, condition (8) holds.

Next two examples answer questions posed by Persi Diaconis and David Freedman.

EXAMPLE 3.14. Let $(\mathcal{X}, d_1)$ and $(\mathcal{Y}, d_2)$ be metric spaces and $\Omega = \mathcal{X} \times \mathcal{Y}$ be equipped with any one of the usual distances

$$d(\omega, \omega^*) = \sqrt{d_1(x, x^*)^2 + d_2(y, y^*)^2}, \qquad d(\omega, \omega^*) = d_1(x, x^*) \vee d_2(y, y^*),$$
$$d(\omega, \omega^*) = d_1(x, x^*) + d_2(y, y^*) \qquad \text{where } \omega = (x, y) \text{ and } \omega^* = (x^*, y^*).$$

By Corollary 3.7, under any such $d$, the balls in $\Omega$ are TIP. Let $D_1 \subset \mathcal{X}$ and $D_2 \subset \mathcal{Y}$ be countable subsets and $(x_1, y_1), (x_2, y_2), \ldots$ any enumeration of the points of $D_1 \times D_2$. Suppose $\{f > 0\} = \bigcup_n H_n$, where $H_n \in \mathcal{F}$ is an open ball centered at $(x_n, y_n)$. For some $k$, the ball $H_k$ is centered at $(x_1, y_2)$. Then, $H_1$ communicates with $H_k$ and $H_k$ communicates with $H_2$, and we let $j_1 = 1$, $j_2 = k$ and $j_3 = 2$. Next, for some $m$, the ball $H_m$ is centered at $(x_2, y_3)$. Then, $H_2$ communicates with $H_m$ and $H_m$ communicates with $H_3$, and we let $j_4 = m$ and $j_5 = 3$. Arguing in this way, $\{f > 0\}$ can be written as $\{f > 0\} = \bigcup_n H_{j_n}$ where $H_{j_n}$ communicates with $H_{j_{n+1}}$ for all $n$. Since each $H_{j_n}$ is TIP, Corollary 3.9 implies $\mathcal{D} = \mathcal{N}$.

EXAMPLE 3.15. In the notation of Example 3.14, suppose $\{f > 0\} = (\bigcup_n H_n)^c$. Let $I_n = \{x : (x, y) \in H_n \text{ for some } y\}$ be the projection of $H_n$ on $\mathcal{X}$. Since $H_n$ is open, $I_n$ is open as well. Suppose also that $\sum_n \mu(I_n) < \mu(\mathcal{X})$. Then, Corollary 3.7 yields $\mathcal{D} = \mathcal{N}$. In fact, $\{f = 0\} \subset (\bigcup_n I_n) \times \mathcal{Y}$ and $\mu(\bigcup_n I_n) \leq \sum_n \mu(I_n) < \mu(\mathcal{X})$. Letting $U_0 = \mathcal{X} - (\bigcup_n I_n)$, thus, one obtains $\{f > 0\} \supset U_0 \times \mathcal{Y}$ and $P(X \in U_0) > 0$.

Let us turn now to $\mathcal{D} \neq \mathcal{N}$. Generally, the complement of a TIP set need not be TIP. One consequence is that, in spite of Corollary 3.8, the intersection of a decreasing sequence of TIP sets need not be TIP.

EXAMPLE 3.16. Let $\mathcal{X} = \mathcal{Y} = (0, 1)$, $\mu = \nu =$ Lebesgue measure, $F = (0, \frac{1}{2}) \times (0, \frac{1}{2})$ and $G_n = (\frac{1}{2} - \frac{1}{n}, 1) \times (\frac{1}{2}, 1)$. Since $F$ and $G_n$ are TIP and $F$ communicates with $G_n$, then $H_n = F \cup G_n$ is TIP. Further, $H_n$ is a decreasing sequence of sets. However,

$$H = \bigcap_n H_n = F \cup ([\tfrac{1}{2}, 1) \times (\tfrac{1}{2}, 1))$$

is not TIP. In fact, $0 < Q_H(F) < 1$, $Q_H(H) = 1$ and

$$F = ((0, \tfrac{1}{2}) \times (0, 1)) \cap H = ((0, 1) \times (0, \tfrac{1}{2})) \cap H.$$

Finally, we exhibit a situation where $\mathcal{D} \neq \mathcal{N}$ though $P$ is absolutely continuous (with respect to Lebesgue measure) and has full topological support.

EXAMPLE 3.17. Let $\mathcal{X} = \mathcal{Y} = (0,1)$ and $\mu = \nu =$ Lebesgue measure. Suppose $\{f > 0\} = \{(x,y) : x,y \in I \text{ or } x,y \in (0,1) - I\}$, where $I$ is a Borel subset of $(0,1)$ satisfying

$$0 < \mu(I \cap J) < \mu(J) \quad \text{for each nonempty open set } J \subset (0,1).$$

Since $0 < P(I \times I) < 1$ and

$$I \times I = (I \times (0,1)) \cap \{f > 0\} = ((0,1) \times I) \cap \{f > 0\},$$

then $\mathcal{D} \neq \mathcal{N}$. Moreover, $P(J_1 \times J_2) \geq P(I \cap J_1 \times I \cap J_2) > 0$ whenever $J_1, J_2 \subset (0,1)$ are nonempty open sets, since $\mu(I \cap J_i) > 0$ for $i = 1, 2$. Thus, $P$ has full topological support.

**4. Two-component Gibbs sampler.** The Gibbs sampler, also known as Glauber dynamics, plays an important role in scientific computing. A detailed treatment can be found in various papers or textbooks; see, for example, [5, 7, 8, 10] and references therein. In this section, the Gibbs-chain is shown to meet the SLLN (Gibbs-admissibility of $P$) if and only if condition (6) holds. Moreover, under mild conditions ($\mathcal{F}$ countably generated and $P \ll \mu \times \nu$), condition (6) is also equivalent to ergodicity of the Gibbs-chain on a certain subset $S_0 \subset \Omega$.

In order to define the Gibbs sampler, $Y$ is assumed to admit a regular version of the conditional distribution given $X$, say $\alpha = \{\alpha(x) : x \in \mathcal{X}\}$. Thus: (i) $\alpha(x)$ is a probability measure on $\mathcal{V}$ for $x \in \mathcal{X}$; (ii) $x \mapsto \alpha(x)(V)$ is $\mathcal{U}$-measurable for $V \in \mathcal{V}$; (iii) $P(U \times V) = \int_U \alpha(x)(V) P \circ X^{-1}(dx)$ for $U \in \mathcal{U}$ and $V \in \mathcal{V}$. Similarly, $X$ is supposed to admit a regular version of the conditional distribution given $Y$, say $\beta = \{\beta(y) : y \in \mathcal{Y}\}$.

The Gibbs-chain $(X_n, Y_n)_{n \geq 0}$ has been informally described in Section 2.3. Formally, $(X_n, Y_n)$ is the homogeneous Markov chain with state space $(\Omega, \mathcal{F})$ and transition kernel

$$K(\omega, U \times V) = K((x,y), U \times V) = \int_V \beta(b)(U) \alpha(x)(db)$$

where $U \in \mathcal{U}, V \in \mathcal{V}$ and $\omega = (x,y) \in \Omega$.

Note that $P$ is a stationary distribution for the chain $(X_n, Y_n)$. Denote $\mathbb{P}$ the law of $(X_n, Y_n)$ when $(X_0, Y_0) \sim P$, and $P_\omega$ the law of $(X_n, Y_n)$ given that $(X_0, Y_0) = \omega$.

Any distributional requirement of $(X_n, Y_n)$ (such as SLLN, CLT, ergodicity, rate of convergence) depends only on the choice of the conditional distributions $\alpha$ and $\beta$. At least if $\mathcal{U}$ and $\mathcal{V}$ are countably generated, however, $\alpha$ and $\beta$ are determined by $P$ up to null sets. Thus, one can try to characterize properties of $(X_n, Y_n)$ via properties of $P$. Here, we first focus on the SLLN and then on ergodicity. Recall that, in Section 2.3, $P$ has been called *Gibbs-admissible* if

$$m_n(\phi) \to \int \phi \, dP \qquad \mathbb{P}\text{-a.s., for all } \phi \in L_1(P)$$

where $m_n(\phi) = \frac{1}{n} \sum_{i=0}^{n-1} \phi(X_i, Y_i)$.

We need the following result.

THEOREM 4.1 [6]. *Given a bounded measurable function $\phi : \Omega \to \mathbb{R}$, define $\phi_0 = \phi$, $\mathcal{G}_n = \sigma(X)$ or $\mathcal{G}_n = \sigma(Y)$ as $n$ is even or odd, and $\phi_n = E(\phi_{n-1} \mid \mathcal{G}_n)$. Then,*

$$E(\phi(X_n, Y_n) \mid (X_0, Y_0) = \omega) = \phi_{2n}(\omega) \qquad \text{for all } n \text{ and } P\text{-almost all } \omega.$$

The previous Theorem 4.1 is a version of Theorem 4.1 of [6]. In the latter paper, the authors focus on densities so that $P$ is assumed absolutely continuous with respect to a product measure. However, such an assumption can be dropped, as it is easily seen from the proof given in [6].

In view of Theorem 4.1 and the Burkholder–Chow result mentioned in Section 2.1, Gibbs-admissibility and $\mathcal{D} = \mathcal{N}$ look like very close conditions. In fact, they are exactly the same thing.

THEOREM 4.2.

$P$ *is Gibbs-admissible* $\Leftrightarrow$ $\mathcal{D} = \mathcal{N}$ $\Leftrightarrow$ *condition* (6) *holds.*

PROOF. The equivalence between (6) and $\mathcal{D} = \mathcal{N}$ has been already proved in Corollary 3.5.

Suppose that $P$ is Gibbs-admissible. In order to check (6), fix $U \in \mathcal{U}$ and $V \in \mathcal{V}$ with $P(U \times V) = P(U^c \times V^c) = 0$. Then, $P(U \times V^c) > 0$ or $P(U^c \times V) > 0$, say $P(U \times V^c) > 0$. Let

$$M_1 = \{y \in V^c : \beta(y)(U^c) > 0\}, \qquad U_1 = \{x \in U : \alpha(x)(V) = \alpha(x)(M_1) = 0\}.$$

Then, $P(U^c \times V^c) = 0$ yields $P(Y \in M_1) = 0$, and $P(Y \in M_1) = 0$ together with $P(U \times V) = 0$ imply $P(X \in U - U_1) = 0$. By induction, for each $j \geq 2$, define

$$M_j = \{y \in V^c : \beta(y)(U_{j-1}^c) > 0\}, \qquad U_j = \{x \in U_{j-1} : \alpha(x)(M_j) = 0\}$$

and verify that $P(Y \in M_j) = 0 = P(X \in U_{j-1} - U_j)$. Define further $U_\infty = \bigcap_j U_j$ and note that $P(X \in U_\infty) = P(X \in U)$. Fix $\omega = (x, y) \in U_\infty \times V^c$. Given $j$, since $\alpha(x)(M_{j+1}) = 0$ and $\beta(b)(U_j) = 1$ for each $b \in V^c - M_{j+1}$, the transition kernel $K$ satisfies

$$K(\omega, U_j \times V^c) = \int_{V^c} \beta(b)(U_j) \alpha(x)(db)$$

$$= \int_{V^c - M_{j+1}} \beta(b)(U_j) \alpha(x)(db) = \alpha(x)(V^c) = 1.$$

Thus, $K(\omega, U_\infty \times V^c) = 1$ for all $\omega \in U_\infty \times V^c$ and this implies

$$P_\omega(X_n \in U_\infty, Y_n \in V^c \text{ for all } n) = 1 \qquad \text{for all } \omega \in U_\infty \times V^c.$$

Next, by Gibbs-admissibility of $P$ (with $\phi = I_{U \times V^c}$), there is a set $S \in \mathcal{F}$ with $P(S) = 1$ and

$$\lim_n m_n(I_{U \times V^c}) = P(U \times V^c) \qquad P_\omega\text{-a.s., for all } \omega \in S.$$

Since $P(S \cap (U_\infty \times V^c)) = P(U \times V^c) > 0$, there is a point $\omega_0 \in S \cap (U_\infty \times V^c)$. For such an $\omega_0$, one obtains

$$P(U \times V^c) = \lim_n m_n(I_{U \times V^c}) = 1 \qquad P_{\omega_0}\text{-a.s.}$$

Therefore $P(Y \in V) = 0$, that is, condition (6) holds.

Finally, suppose $\mathcal{D} = \mathcal{N}$. By the ergodic theorem, since $(X_n, Y_n)$ is stationary under $\mathbb{P}$, for $P$ to be Gibbs-admissible it is enough that $\mathbb{P}$ be 0–1-valued on the shift–invariant sub-$\sigma$-field of $\mathcal{F}^\infty$. Let $h$ be a bounded harmonic function, that is, $h : \Omega \to \mathbb{R}$ is bounded measurable and $h(\omega) = \int h(t) K(\omega, dt)$ for all $\omega \in \Omega$. Because of Theorem 4.1 and $h$ harmonic,

$$h(\omega) = E(h(X_n, Y_n) \mid (X_0, Y_0) = \omega) = h_{2n}(\omega) \qquad \text{for all } n \text{ and } P\text{-almost all } \omega.$$

By the Burkholder–Chow result (Section 2.1) and $\mathcal{D} = \mathcal{N}$, one also obtains

$$h_{2n} \to E(h \mid \mathcal{D}) = \int h \, dP \qquad P\text{-a.s.}$$

Hence, $h(\omega) = \int h \, dP$ for $P$-almost all $\omega$. Let $H \in \mathcal{F}^\infty$ be such that $H = \theta^{-1} H$, where $\theta$ is the shift transformation on $\Omega^\infty$. Then,

$$h(\omega) = P_\omega(H)$$

is a bounded harmonic function satisfying $I_H = \lim_n h(X_n, Y_n)$, $\mathbb{P}$-a.s.; see, for example, Theorem 17.1.3 of [9]. Since $(X_n, Y_n)$ is stationary under $\mathbb{P}$, then $I_H = h(X_0, Y_0)$, $\mathbb{P}$-a.s. Hence, $P(h = 0) = 1$ or $P(h = 1) = 1$, which implies $\mathbb{P}(H) = \int h \, dP \in \{0, 1\}$. This concludes the proof. $\square$

Theorem 4.2 is potentially useful in real problems as well, since it singles out those $P$ such that Gibbs sampling makes sense; see Section 2.3.

The next example is motivated by, in Gibbs sampling applications, the available information typically consisting of the conditionals $\alpha$ and $\beta$.

EXAMPLE 4.3. Condition (6) can be stated in terms of the conditional distribution $\alpha = \{\alpha(x) : x \in \mathcal{X}\}$ of $Y$ given $X$. Let $V \in \mathcal{V}$. If $P(U \times V) = P(U^c \times V^c) = 0$ for some $U \in \mathcal{U}$, then $\alpha(X)(V) = I_{U^c}(X) \in \{0, 1\}$ a.s. Conversely, $\alpha(X)(V) \in \{0, 1\}$ a.s. implies $P(U \times V) = P(U^c \times V^c) = 0$ with $U = \{x : \alpha(x)(V) = 0\}$. It follows that (6) is equivalent to

(11) $\quad \alpha(X)(V) \in \{0, 1\}$ a.s. $\implies$ $\alpha(X)(V) = 0$ a.s. or $\alpha(X)(V) = 1$ a.s.

for all $V \in \mathcal{V}$. Thus, whether or not $P$ is Gibbs-admissible depends only on the "supports" of the probability laws $\alpha(x)$, $x \in \mathcal{X}$. Condition (11) also suggests various sufficient criteria. Let

$$\mathcal{V}_0 = \{V \in \mathcal{V} : 0 < P(\alpha(X)(V) = 1) < 1\}.$$

Condition (11) trivially holds whenever $V \notin \mathcal{V}_0$. Hence, a first (obvious) sufficient condition for Gibbs-admissibility of $P$ is

(12) $\qquad \alpha(X)(V) > 0 \qquad$ a.s. for each $V \in \mathcal{V}_0$.

A second condition is the following. Let $\mathcal{X}$ be a metric space and $\mathcal{U}$ the Borel $\sigma$-field. Then, $P$ is Gibbs-admissible provided there is a set $T \in \mathcal{U}$ satisfying:

(i) $P(X \in T) = 1$ and $P(X \in U) > 0$ for all open $U \subset \mathcal{X}$ such that $T \cap U \neq \varnothing$;
(ii) for each $V \in \mathcal{V}_0$, the map $x \mapsto \alpha(x)(V)$ is continuous on $T$;
(iii) for each $V \in \mathcal{V}_0$, there is $x \in T$ such that $0 < \alpha(x)(V) < 1$.

Fix in fact $V \in \mathcal{V}_0$. Then, (ii)–(iii) imply $\{x \in T : 0 < \alpha(x)(V) < 1\} = T \cap U \neq \varnothing$ for some open $U \subset \mathcal{X}$. Thus, $P(0 < \alpha(X)(V) < 1) = P(X \in U) > 0$ by (i). Therefore, condition (11) holds. Note that, in view of (ii), condition (iii) is certainly true if $T$ is connected. Similarly, (iii) holds if, for each $V \in \mathcal{V}_0$, the map $x \mapsto \alpha(x)(V)$ is not constant on some connected component of $T$.

In applications, it is useful that $(X_n, Y_n)$ is ergodic on some *known* set $S \in \mathcal{F}$. By ergodicity on $S$, we mean $S \in \mathcal{F}$ and

$$P(S) = 1, \qquad K(\omega, S) = 1 \quad \text{and} \quad \|K^n(\omega, \cdot) - P\| \to 0 \qquad \text{for all } \omega \in S,$$

where $\|\cdot\|$ is total variation norm and $K^n$ the $n$th iterate of $K$. If $(X_n, Y_n)$ is ergodic on $S$, for each $\omega \in S$ one obtains

$$m_n(\phi) \to \int \phi \, dP \qquad P_\omega\text{-a.s., for all } \phi \in L_1(P).$$

Thus, ergodicity on some $S$ implies Gibbs-admissibility of $P$. We now seek conditions for the converse to be true.

To this end, an intriguing choice of $S$ is

$$S_0 = \{\omega \in \Omega : K(\omega, \cdot) \ll P\}.$$

A simple condition for $S_0 \in \mathcal{F}$ is $\mathcal{F}$ countably generated.

THEOREM 4.4. *If $\mathcal{F}$ is countably generated, condition* (6) *holds and $P(S_0) = 1$, then $(X_n, Y_n)$ is ergodic on $S_0$.*

PROOF. Since $P(S_0) = 1$, the definition of $S_0$ gives $K(\omega, S_0) = 1$ for all $\omega \in S_0$. Let $P_0$ and $K_0(\omega, \cdot)$ be the restrictions of $P$ and $K(\omega, \cdot)$ to $\mathcal{F}_0$, where $\omega \in S_0$ and $\mathcal{F}_0 = \{F \cap S_0 : F \in \mathcal{F}\}$. Then, $(X_n, Y_n)$ can be seen as a Markov chain with state space $(S_0, \mathcal{F}_0)$, transition kernel $K_0$ and stationary distribution $P_0$. Also,

$$\|K^n(\omega, \cdot) - P\| = \|K_0^n(\omega, \cdot) - P_0\| \qquad \text{for all } \omega \in S_0.$$

By standard arguments on Markov chains, thus, it is enough to prove that $K_0$ is aperiodic and every bounded harmonic function (with respect to $K_0$) is constant on $S_0$. Let $h_0$ be one such function, that is, $h_0 : S_0 \to \mathbb{R}$ is bounded measurable and $h_0(\omega) = \int h_0(t) K_0(\omega, dt)$ for all $\omega \in S_0$. Define $h = h_0$ on $S_0$ and $h = 0$ on $S_0^c$. Then, $h : \Omega \to \mathbb{R}$ is bounded measurable and $h(\omega) = \int h(t) K(\omega, dt)$ for $P$-almost all $\omega$. Letting $A = \{\omega \in \Omega : h(\omega) = \int h \, dP\}$ and arguing as in the proof of Theorem 4.2, condition (6) implies $P(A) = 1$. Thus, $K(\omega, A) = 1$ for each $\omega \in S_0$, so that

$$h_0(\omega) = \int h_0(t) K_0(\omega, dt) = \int_A h(t) K(\omega, dt) = \int h \, dP \qquad \text{for all } \omega \in S_0.$$

It remains to prove aperiodicity of $K_0$. Toward a contradiction, suppose there are $d \geq 2$ nonempty disjoint sets $F_1, \ldots, F_d \in \mathcal{F}_0$ such that

$$K_0(\omega, F_{i+1}) = 1 \qquad \text{for all } \omega \in F_i \text{ and } i = 1, \ldots, d, \text{ where } F_{d+1} = F_1.$$

If $P(F_1) = 1$, then $K_0(\omega, F_1) = 1$ for all $\omega \in S_0$, contrary to $K_0(\omega, F_1) = 0$ for $\omega \in F_1$. Hence, $P(F_1) < 1$. Applying Theorem 4.1 to $\phi = I_{F_1}$, one obtains

$$K^{nd}(\omega, F_1) = \phi_{2nd}(\omega) \qquad \text{for all } n \text{ and } P\text{-almost all } \omega.$$

Hence, the Burkholder–Chow result (Section 2.1) and condition (6) yield

$$K^{nd}(\cdot, F_1) = \phi_{2nd} \to E(\phi \mid \mathcal{D}) = \int \phi \, dP = P(F_1) \qquad \text{a.s.}$$

Since $\lim_n K^{nd}(\omega, F_1) = 1 \neq P(F_1)$ for all $\omega \in F_1$, it follows that $P(F_1) = 0$. But this is a contradiction, since $P(F_1) = 0$ implies $K_0(\omega, F_1) = 0$ for all $\omega \in S_0$. Thus, $K_0$ is aperiodic. $\square$

By Theorem 4.4, Gibbs-admissibility implies ergodicity on $S_0$ whenever $P(S_0) = 1$ (and $\mathcal{F}$ is countably generated). In turn, for $P(S_0) = 1$, it is enough that $P \ll \mu \times \nu$.

THEOREM 4.5. *If $\mathcal{F}$ is countably generated and $P \ll \mu \times \nu$, then $(X_n, Y_n)$ is ergodic on $S_0$ if and only if condition (6) holds, that is, if and only if $P$ is Gibbs-admissible.*

PROOF. It suffices to prove $P(S_0) = 1$. Let $f$ be a version of the density of $P$ with respect to $\mu \times \nu$ and

$$f_1(x) = \int f(x,y)\nu(dy), \qquad f_2(y) = \int f(x,y)\mu(dx).$$

Define $D_1 = \{x : 0 < f_1(x) < \infty\}$, $D_2 = \{y : 0 < f_2(y) < \infty\}$ and

$$D = \{x \in D_1 : \alpha(x)(D_2) = 1\}.$$

Since $P(X \in D_1) = P(Y \in D_2) = 1$, then $P(X \in D) = 1$. Fix $\omega = (x,y) \in D \times \mathcal{Y}$. Then, $\alpha(x)$ has density $f(x,\cdot)/f_1(x)$ with respect to $\nu$. Also, $\beta(b)$ has density $f(\cdot,b)/f_2(b)$ with respect to $\mu$ for each $b \in D_2$, and $\alpha(x)(D_2) = 1$. Hence, for $C \in \mathcal{F}$,

$$K(\omega, C) = \int \int I_C(a,b)\beta(b)(da)\alpha(x)(db)$$
$$= \int \int I_C(a,b)\frac{f(a,b)}{f_2(b)}\mu(da)\frac{f(x,b)}{f_1(x)}\nu(db)$$
$$= \frac{1}{f_1(x)}\int \int I_C(a,b)\frac{f(x,b)}{f_2(b)}f(a,b)\mu(da)\nu(db)$$
$$= \frac{1}{f_1(x)}\int I_C(a,b)\frac{f(x,b)}{f_2(b)}P(d(a,b)).$$

Therefore, $K(\omega, \cdot) \ll P$, so that $P(S_0) \geq P(D \times \mathcal{Y}) = 1$. □

EXAMPLE 4.3 (Continued). Suppose $\mathcal{F}$ is countably generated and $P \ll \mu \times \nu$. By Theorem 4.5, $(X_n, Y_n)$ is ergodic on $S_0$ if and only if $\alpha$ meets condition (11). Let $f$ and $f_1$ be as in the proof of Theorem 4.5 and

$$I_x = \{y : f(x,y) > 0\}.$$

Then, condition (11) holds provided

$$\nu(I_x \cap I_z) > 0 \qquad \text{whenever } 0 < f_1(x), f_1(z) < \infty.$$

Suppose in fact $\alpha(x)(V) = 1$ for some $x$ with $0 < f_1(x) < \infty$ and $V \in \mathcal{V}_0$. For every $z$ satisfying $0 < f_1(z) < \infty$, one obtains $\alpha(z)(I_x) > 0$ [since $\nu(I_x \cap I_z) > 0$] and this implies $\alpha(z)(V) > 0$. Hence, condition (12) holds. Next, suppose $\mathcal{X}$ is a metric space, $\mathcal{U}$ the Borel $\sigma$-field and $\mu$ has full topological support. Then, another sufficient condition for (11) is:

(j) $f(x,y) \leq h(y)$ for all $(x,y)$ and some $\nu$-integrable function $h$;
(jj) $x \mapsto f(x,y)$ is continuous for fixed $y \in \mathcal{Y}$;
(jjj) for each $V \in \mathcal{V}_0$, there is $x \in \mathcal{X}$ with $\nu(I_x \cap V) > 0$ and $\nu(I_x \cap V^c) > 0$.

Under (j)–(jj), in fact, $f_1$ is a real continuous function and $x \mapsto \alpha(x)(V)$ is continuous on the set $\{f_1 > 0\}$ for all $V \in \mathcal{V}_0$. Also, since $\{f_1 > 0\}$ is open and $\mu$ has full topological support, $P(X \in U) > 0$ whenever $U$ is open and $U \cap \{f_1 > 0\} \neq \varnothing$. Hence, conditions (i)–(iii) (mentioned in the first part of this example) are satisfied with $T = \{f_1 > 0\}$. Recall that (jjj) holds if $\{f_1 > 0\}$ is connected (as well as in some other situations).

We close the paper with two remarks.

REMARK 4.6 (Uniform and geometric ergodicity). Let $f$, $f_1$, $f_2$ and $D$ be as in the proof of Theorem 4.5, where $\mathcal{F}$ is countably generated and $P \ll \mu \times \nu$. Suppose that

$$sI_{U \times V} \leq f \quad \text{and} \quad f_1 \leq tI_U$$

for some constants $s, t > 0$ and $U \in \mathcal{U}$, $V \in \mathcal{V}$ with $P(U \times V) > 0$.

Then, $(X_n, Y_n)$ is ergodic on $S_0$. In addition, $(X_n, Y_n)$ is *uniformly ergodic* on $D \times \mathcal{Y}$, in the sense that

$$\sup_{\omega \in D \times \mathcal{Y}} \|K^n(\omega, \cdot) - P\| \leq qr^n$$

for some constants $q > 0$ and $r < 1$. Also, the convergence rate $r$ can be taken such that $r \leq 1 - (s/t)\nu(V)$.

In fact, $f_1 = 0$ on $U^c$ implies $f = 0$ on $U^c \times \mathcal{Y}$, $\mu \times \nu$-a.e. Thus, (6) holds by Corollaries 3.5 and 3.7, and $(X_n, Y_n)$ is ergodic on $S_0$ by Theorem 4.5. Since $\mu(U) > 0$, $\nu(V) > 0$ and

$$s\mu(U)\nu(V) \leq \int_{U \times V} f \, d(\mu \times \nu) = P(U \times V),$$

then $0 < \nu(V) < \infty$, and one can define the probability measure

$$\gamma(C) = \frac{1}{\nu(V)} \int I_C(a,b) I_V(b) \frac{1}{f_2(b)} P(d(a,b)), \qquad C \in \mathcal{F}.$$

Since $f_1 = 0$ on $U^c$, then $D \subset \{f_1 > 0\} \subset U$. Therefore, for each $\omega = (x, y) \in D \times \mathcal{Y}$, one obtains

$$K(\omega, C) = \frac{1}{f_1(x)} \int I_C(a,b) \frac{f(x,b)}{f_2(b)} P(d(a,b))$$

$$\geq \frac{s}{t} \int I_C(a,b) I_V(b) \frac{1}{f_2(b)} P(d(a,b)) = \frac{s}{t} \nu(V) \gamma(C), \qquad C \in \mathcal{F}.$$

Thus, $D \times \mathcal{Y}$ is a *small set* such that $P(D \times \mathcal{Y}) = 1$, and this implies uniform ergodicity of $(X_n, Y_n)$ on $D \times \mathcal{Y}$; see pages 1714–1715 and Proposition 2 of [10].

The previous assumptions can be adapted to obtain *geometric ergodicity*, in the sense that $(X_n, Y_n)$ is ergodic on $S_0$ and $\|K^n(\omega, \cdot) - P\| \leq q(\omega) r^n$ for $P$-almost all $\omega$, where $r \in (0,1)$ is a constant and $q$ a function in $L_1(P)$. As an example (we omit calculations), $(X_n, Y_n)$ is geometrically ergodic whenever

$f$ and $f_1$ are bounded, $\quad f \geq s$ on $U \times V$, $\quad f = 0$ on $U^c \times V^c$,

for some $s > 0$, $U \in \mathcal{U}$, $V \in \mathcal{V}$ such that

$$P(U \times V) > 0 \quad \text{and} \quad \sup_{\omega \in U^c \times V} f(\omega) < s \frac{\mu(U)}{\mu(U^c)}.$$

Note that, for the above conditions to apply, $\mu$ must be a finite measure. Even if long to be stated, such conditions can be useful. They apply, for instance, when $(\Omega, \mathcal{F})$ is the Borel unit square, $\mu = \nu =$ one-dimensional Lebesgue measure, and $P$ uniform on the lower or upper half triangle.

REMARK 4.7 (The $k$-component case). This paper has been thought and written for the two-component case, but its contents extend to the $k$-component case, with $k \geq 2$ any integer. In particular, Theorems 4.2, 4.4 and 4.5 can be adapted to the $k$-component Gibbs sampler. We just mention that, in general, the involved sub-$\sigma$-fields are $\mathcal{A}_i = \sigma(Z_1, \ldots, Z_{i-1}, Z_{i+1}, \ldots, Z_k)$, $i = 1, \ldots, k$, where $Z_i$ denotes the $i$th coordinate projection on the product of some $k$ measurable spaces.

**Acknowledgments.** We are grateful to Persi Diaconis for suggesting the problem and for encouraging and helping us.

P. Berti
Dipartimento di Matematica Pura
 ed Applicata "G. Vitali"
Universita' di Modena e Reggio-Emilia
via Campi 213/B
41100 Modena
Italy
E-mail: patrizia.berti@unimore.it

L. Pratelli
Accademia Navale
viale Italia 72
57100 Livorno
Italy
E-mail: pratel@mail.dm.unipi.it

P. Rigo
Dipartimento di Economia Politica
 e Metodi Quantitativi
Universita' di Pavia
via S. Felice 5
27100 Pavia
Italy
E-mail: prigo@eco.unipv.it